\documentclass[a4paper,12pt]{amsart}
\usepackage{amsfonts}
\usepackage{amssymb}
\usepackage{ifthen}
\usepackage{graphicx}
\usepackage{color}
\nonstopmode \numberwithin{equation}{section}
\newtheorem{thm}{Theorem}%[section]
\newtheorem{lem}{Lemma}%[section]
\newtheorem{cor}{Corollary}%[section]
\newtheorem{prop}{Proposition}%[section]

\newtheorem{conj}{Conjecture}

\theoremstyle{definition}
\newtheorem{defn}{Definition}%[section]
\newtheorem{example}{Example}%[section]

\newtheorem{ques}{Question}
\newtheorem{rem}{Remark}

\newtheorem{rems}{Remarks}

\newcounter {own}
\def\theown {\thesection  .\arabic{own}}

\newenvironment{pf}[1][]{%
 \vskip 3mm
 \noindent
 \ifthenelse{\equal{#1}{}}%
  {{\slshape Proof. }}%
  {{\slshape #1.} }%
 }%
{\qed\bigskip}

\newcounter{alphabet}
\newcounter{tmp}

\makeatletter
\newcommand{\Ref}[1]{\@ifundefined{r@#1}{}{\setcounter{tmp}{\ref{#1}}\Alph{tmp}}}
\makeatother
\newenvironment{Lem}[1][]{\refstepcounter{alphabet}%
\bigskip%
\noindent%
{\bf Lemma \Alph{alphabet}}%
{\bf .} \itshape}{\vskip 8pt}

\newcommand{\N}{{\mathbb N}}
\newcommand{\C}{{\mathbb C}}
\newcommand{\D}{{\mathbb D}}

\newcommand{\T}{{\mathbb{T}}}

\newcommand{\Aut}{{\operatorname{Aut}}}

\def\be{\begin{equation}}
\def\ee{\end{equation}}

\newcommand{\bee}{\begin{enumerate}}
\newcommand{\eee}{\end{enumerate}}

\newcommand{\blem}{\begin{lem}}
\newcommand{\elem}{\end{lem}}
\newcommand{\bthm}{\begin{thm}}
\newcommand{\ethm}{\end{thm}}
\newcommand{\bcor}{\begin{cor}}
\newcommand{\ecor}{\end{cor}}
\newcommand{\beg}{\begin{example}}
\newcommand{\eeg}{\end{example}}
\newcommand{\begs}{\begin{examples}}
\newcommand{\eegs}{\end{examples}}
\newcommand{\bdefn}{\begin{defn}}
\newcommand{\edefn}{\end{defn}}
\newcommand{\bprob}{\begin{prob}}
\newcommand{\eprob}{\end{prob}}
\newcommand{\bei}{\begin{itemize}}
\newcommand{\eei}{\end{itemize}}
\newcommand{\bqn}{\begin{ques}}
\newcommand{\eqn}{\end{ques}}
\newcommand{\bcon}{\begin{conj}}
\newcommand{\econ}{\end{conj}}
\newcommand{\bcons}{\begin{conjs}}
\newcommand{\econs}{\end{conjs}}
\newcommand{\bprop}{\begin{prop}}
\newcommand{\eprop}{\end{prop}}
\newcommand{\brem}{\begin{rem}}
\newcommand{\erem}{\end{rem}}
\newcommand{\brems}{\begin{rems}}
\newcommand{\erems}{\end{rems}}
\newcommand{\bo}{\begin{obser}}
\newcommand{\eo}{\end{obser}}
\newcommand{\bos}{\begin{obsers}}
\newcommand{\eos}{\end{obsers}}
\newcommand{\bpf}{\begin{pf}}
\newcommand{\epf}{\end{pf}}
\newcommand{\ba}{\begin{array}}
\newcommand{\ea}{\end{array}}
\newcommand{\beq}{\begin{eqnarray}}
\newcommand{\beqq}{\begin{eqnarray*}}
\newcommand{\eeq}{\end{eqnarray}}
\newcommand{\eeqq}{\end{eqnarray*}}

\newcommand{\Ra}{\Rightarrow}

\newcommand{\ds}{\displaystyle}

\newcounter{minutes}
\divide\time by 60
\newcounter{hours}
\multiply\time by 60 \addtocounter{minutes}{-\time}
\begin{document}
\bibliographystyle{amsplain}
\title[Composition operators on the discrete Hardy space on homogenous trees]
{Composition operators on the discrete analogue of generalized Hardy space on homogenous trees}

%%%=========================================================================
%%\thanks{%$^\dagger$
%File:~\jobname .tex,
%          printed: \number\day-\number\month-\number\year,
%          \thehours.\ifnum\theminutes<10{0}\fi\theminutes}
%%%=========================================================================
%
%\author{S. Ponnusamy $^\dagger $}

\author{Perumal Muthukumar}
\address{P. Muthukumar,
Indian Statistical Institute (ISI), Chennai Centre,
SETS (Society for Electronic Transactions and Security),
MGR Knowledge City, CIT Campus, Taramani,
Chennai 600 113, India. }
\email{pmuthumaths@gmail.com}

\author{Saminathan Ponnusamy
}
\address{S. Ponnusamy,
Indian Statistical Institute (ISI), Chennai Centre,
SETS (Society for Electronic Transactions and Security),
MGR Knowledge City, CIT Campus, Taramani,
Chennai 600 113, India.
}
\email{samy@isichennai.res.in, samy@iitm.ac.in}

\subjclass[2000]{Primary: 05C05, 47B33, 47B38; Secondary: 46B50}
\keywords{Composition operators, Rooted homogeneous tree, discrete Hardy spaces.\\
}

%\date{\today

\begin{abstract}
In this paper, we study the basic properties such as boundedness and compactness of composition operators on discrete analogue of generalized Hardy space
defined on a homogeneous rooted tree. Also, we compute  the operator norm of composition operator when inducing symbol is automorphism of a homogenous tree.
\end{abstract}
\thanks{
%%=========================================================================
File:~\jobname .tex,
          printed: \number\day-\number\month-\number\year,
          \thehours.\ifnum\theminutes<10{0}\fi\theminutes
%%=========================================================================
}
\maketitle
\pagestyle{myheadings}
\markboth{P. Muthukumar and S. Ponnusamy}{Composition Operators on the Discrete Hardy Space on Homogenous Trees}

\section{Introduction}\label{MP2Sec1}
Let $\Omega$ be a nonempty set and $X$ be a complex Banach space of complex valued functions defined on $\Omega$.
For a self map $\phi$ of $\Omega$,  the composition operator $C_\phi$ induced by the symbol $\phi$ is defined as
$$ C_\phi(f)=g ~\mbox{ where }~ g(x)=f(\phi(x)) ~\mbox{ for all $x\in \Omega$ and $f\in X$}.
$$
In the classical case, $\Omega$ is the unit disk $\D= \{z \in \C:\, |z|<1\}$ and the choices for $X$ are analytic functions
spaces, eg. the Hardy space $H^p$, the Bergman space $A^p$, the Bloch space $\mathcal B$, etc. The study of composition
operators on various analytic function spaces defined on ${\mathbb D}$ is well known. There are excellent books on
composition operators, see \cite{Cowen:Book,Shapiro:Book,Rksingh:Book} and the references therein. The approach in
the first two books \cite{Cowen:Book,Shapiro:Book} are function theoretic whereas \cite{Rksingh:Book} deals in measure
theoretic point of view. Also, there a number of articles dealing with composition operators on different transform spaces,
see for example \cite{Mu,Mu-Yall,CM}. In this article, $\Omega$ will be a homogeneous rooted tree and $X$ the discrete
analogue of generalized Hardy space introduced in \cite{MP-Tp-spaces}.

In the recent years, there has been a great interest in studying operator theory on discrete structure such as graphs,
in particular on an infinite tree graph \cite{Colonna-MO-5,Colonna-MO-3,Colonna-MO-2,Colonna-MO-1,Colonna-MO-4}.
In \cite{Colonna-CO}, Colonna et al.  studied composition operators on Lipschitz functions on a tree with edge counting
metric to the complex plane with Euclidean metric, which is a discrete analogue of Bloch space, because Bloch space is
also consisting of only Lipschitz functions on the unit disk under Hyperbolic metric to the complex plane with Euclidean metric.
In \cite{MP-Tp-spaces}, the present authors defined discrete analogue ($\mathbb{T}_{p}$) of generalized Hardy spaces on
homogeneous rooted tree. In the same article multiplication operators on $\mathbb{T}_{p}$ spaces are studied.

In this article,  we deal with the study of composition operators on $\mathbb{T}_{p}$ spaces. We refer to  Section \ref{MP2Sec2}
for the definitions of  homogeneous rooted tree and $\mathbb{T}_{p}$ spaces. In Section \ref{MP2Sec4}, we consider the
boundedness of composition operators on $\mathbb{T}_{p}$ spaces and some of its consequences including norm estimates.
In Section \ref{MP2Sec5}, we consider the compactness of composition operators on $\mathbb{T}_{p}$ spaces and
derive equivalent conditions for compactness. Finally, in Section \ref{MP2Sec6}, we present three examples to show the following:
there are self maps of $T$ which do not induce bounded composition operator on $\mathbb{T}_{p}$; there exists a bounded
composition operator on $\mathbb{T}_{p}$ which is not compact; there are unbounded self maps of $T$ which induces compact
composition operators on $\mathbb{T}_{p}$ for the case of $(q+1)$-homogeneous trees with $q\geq 2$.

\section{Preliminaries and Lemmas}\label{MP2Sec2}
Let $G=(V,E)$ be a graph such that $E\subseteq V \times V$, where the elements of the sets $V$ and $E$ are called vertices
and edges of the graph $G$, respectively. We shall not always distinguish between a graph and its vertex set and so,
we may write $x \in G$ (rather than $x \in V$) and by a function defined on a graph, we mean a function defined on its vertices.
Two vertices $x,y \in G$  are said to be \textit{neighbours} (denoted by $x\sim y$) if $(x,y)\in E$. If all the vertices of $G$
have the same number $k$ of neighbours, then the graph is said to be \textit{$k$-homogeneous} or $k$-regular graph.
A \textit{ finite path} is a nonempty subgraph $P=(V,E)$ of the form $V=\{x_0,x_1,\ldots ,x_k\}$ and
$E=\{(x_0,x_1),(x_1,x_2),\ldots,(x_{k-1},x_k)\}$, where $x_i$'s are distinct. In this case, we call $P$ be a path
between $x_0$ and $x_k$. If $P$ is a path between $x_0$ and $x_k ~(k\geq 2)$, then $P$ with an additional
edge $(x_{n},x_{0})$ is called a \textit{cycle}. A nonempty graph $G$ is called \textit{connected} if for any two of its vertices,
there is a path between them. A connected  graph without cycles is called a \textit{tree}. Thus, any two vertices of a
tree are linked by a unique path. The \textit{distance} between any two vertex of a
tree is the number of edges in the unique path connecting them. Sometimes it is convenient to consider one vertex of a tree as
special; such a vertex is then called the root of this tree. A tree $T$ with fixed root $\textsl{o}$ is called a \textit{rooted tree}.
If $G$ is a rooted tree with root $\textsl{o}$, then $|v|$ denotes the distance between the root $\textsl{o}$ and the vertex $v$.
Further, the \textit{parent} (denoted by $v^-$) of a vertex $v$, which is not a root, is the unique vertex $w\in G$ such that
$w\sim v$ and $|w|=|v|-1$. In this case, $v$ is called \textit{child} of $w$.
For basic issues regarding graph theory, one can refer standard texts on this subject.

Throughout the paper, unless otherwise stated explicitly, $T$ denotes a homogeneous rooted tree (hence infinite graph), $\phi$
 denotes a self map of $T$,
$\mathbb{N}=\{1,2,\ldots \}$  and $\mathbb{N}_0=\mathbb{N}\cup \{0\}$.

For $p\in(0,\infty]$, the {\it  Hardy space}
$H^{p}$ consists of all those analytic functions $f:\mathbb{D}\rightarrow\mathbb{C}$ such that
$\|f\|_{p}<\infty$,
where
$$
\|f\|_{p}=\sup_{0\leq r<1}M_{p}(r,f)
$$
and
$$M_{p}(r,f)=
\begin{cases}
\displaystyle \left(\frac{1}{2\pi}\int_{0}^{2\pi}|f(re^{i\theta})|^{p}\,d\theta\right)^\frac{1}{p}
& \mbox{if } p\in(0,\infty)\\
\displaystyle\sup_{|z|=r}|f(z)| &\mbox{if } p=\infty.
\end{cases}
$$
The {\it generalized Hardy space} $H^{p}_{g}$ is defined similarly, upon replacing analytic functions by measurable functions.

As in \cite{MP-Tp-spaces}, in a $(q+1)$-homogeneous tree $T$ rooted at $\textsl{o}$, we define
$$\|f\|_{p}:= \sup\limits_{n\in \mathbb{N}_{0}} M_{p}(n,f),
$$
where $M_{p}(0,f):= |f(\textsl{o})|$ and for every $n\in \mathbb{N}$,
$$%\be\label{MP2eq2}
M_{p}(n,f):=
\left\{
\begin{array}{ll}
\ds \left (\frac{1}{(q+1)q^{n-1}}\sum\limits_{|v|=n}|f(v)|^{p} \right )^{\frac{1}{p}} & \mbox{if } p\in(0,\infty) \\
\max\limits_{|v|=n } |f(v)| & \mbox{if } p=\infty .
\end{array}
\right.
$$
The discrete analogue of the generalized Hardy space, denoted by $\mathbb{T}_{q,p}$, is then defined by
$$\mathbb{T}_{q,p}:=\{f\colon T \to\mathbb{C} \, \big |\, \|f\|_{p}<\infty\}
$$
for every $p\in(0,\infty]$. For the sake of simplicity, we shall write $\mathbb{T}_{q,p}$  as $\mathbb{T}_{p}$.
Throughout the discussion,  $\|.\|$ denotes $\|.\|_p$  in $\mathbb{T}_{p}$ spaces. The following results proved
by the present authors \cite{MP-Tp-spaces} are needed for our present investigation.

\begin{Lem}\label{thm:banachp}
For $1\leq p\leq\infty$, $\|.\|_{p}$ induces a Banach space structure on the space $\mathbb{T}_{p}$.
\end{Lem}

\begin{Lem} \emph{(Growth Estimate)} \label{lem:bound}
Let $T$ be a $(q+1)$-homogeneous tree rooted at $\textsl{o}$ and $0<p<\infty$. Then, for $v\in T$, we have the following:
If $f\in \mathbb{T}_{p}$, then
$$ |f(v)|\leq \{(q+1)q^{|v|-1}\}^{\frac{1}{p}} \|f\|_p.
$$
\end{Lem}

\begin{Lem}
Norm convergence in $\mathbb{T}_{p}$ implies pointwise convergence. That is,
$$\lim\limits_{n\rightarrow\infty}\|f_n-f\| =0~ \Rightarrow ~ \lim\limits_{n\rightarrow\infty}f_n(v) = f(v)
\mbox{ for each $v\in T$.}
$$
\end{Lem}

%Proof of all these results stated in this section are straight forward from the definitions. For the proof, see \cite{MP-Tp-spaces}.

\section{Bounded Composition Operators} \label{MP2Sec4}

A linear operator $A$ from a normed linear space $X$ to a normed linear space $Y$ is said to be \textit{bounded}
if the operator norm $ \|A\|= \mbox{sup}\{\|Ax\|_Y : \|x\|_X=1\}$ is finite.

Before we proceed to discuss our results, it is appropriate to recall some basic results about bounded composition operators in the classical
setting.  For example (see \cite[ Corollary 3.7]{Cowen:Book}), every analytic self map $\phi$ of $\D$ induces bounded
composition operator $C_\phi$ on $H^p$, $1\leq p<\infty$. Moreover,
\be\label{normHp}
\|C_\phi\|^p\leq \frac{1+|\phi(0)|}{1-|\phi(0)|}.
\ee
It is also known that (see \cite[Theorem 3.8]{Cowen:Book})  equality holds in (\ref{normHp}) for every inner function
of $\D$ (for example, for every automorphism of $\D$). For the case $p=\infty$, it is easy to see that $\|C_\phi\|=1$
for every analytic self map $\phi$ of $\D$.

Now, for our setting, we let $\phi$ be a self map of $(q+1)$-homogeneous rooted tree $T$. For $n\in \mathbb{N}_{0}$ and $w\in T$,
let $N_{\phi}(n,w)$ denote the number of pre-images of $w$ for $\phi$ in $|v|=n$. That is
$N_{\phi}(n,w)$ is the number of elements in $\{\phi^{-1}(w)\}\bigcap \{|v|=n\}$.
For $w\in T$, we define the weight function $W$ as follows:
\be\label{MP2-eq1}
W(w):=
\left\{
\begin{array}{ll}
(q+1)q^{|w|-1} & \mbox{if } w\in T\setminus \{\textsl{o}\} \\
1 & \mbox{if } w=\textsl{o}.
\end{array}
\right.
\ee
Let $|D_n|$  denote the number of vertices with $|v|=n$. Thus,
$$%\be\label{}
|D_n|=\left\{\begin{array}{ll}
(q+1)q^{n-1} & \mbox{if } n\in \mathbb{N} \\
1 & \mbox{if } n=0.
\end{array} \right.
$$
\bthm
Every self map $\phi$ of $T$ induces bounded composition operator on $\T_\infty$ with $\|C_\phi\|=1$.
\ethm
\bpf
For each $f\in\T_\infty$ and every self map $\phi$ of $T$, we have
$$
\|C_\phi(f)\|_\infty=\|f\circ\phi\|_\infty=\sup\limits_{w\in \phi(T)}|f(w)|\leq\|f\|_\infty.
$$
Thus, $C_\phi$ is bounded on $\T_\infty$. It is easy to see that $\|\chi_{\{v\}}\circ\phi\|_\infty=1$ for each
$v\in T$, where $\chi_{\{v\}}$ denotes the characteristic function on $\{v\}$. It gives that $\|C_\phi\|=1$.
\epf

In order to study the boundedness of the composition operators on $\T_p$ for $1\leq p<\infty$, it is convenient to
deal with the case $q=1$ and $q\geq 2$ independently. First, we begin with the case $q=1$.

\bthm\label{MP2-th1}
For every self map $\phi$ of $2$-homogeneous tree $T$, $C_\phi$ is bounded on $\mathbb{T}_{p}$ with $\|C_\phi\|^p \leq 2$, $1\leq p<\infty$.
\ethm
\bpf
By the growth estimate  (Lemma \Ref{lem:bound}) for $2$-homogeneous trees, it follows that
$|f(v)|^p \leq 2\|f\|^p$ for all $v\in T$ and  $f\in \mathbb{T}_{p}$.
So,$$M_p^p(0,C_\phi f)=|f(\phi(\textsl{o}))|^p \leq 2\|f\|^p .$$  Since $|D_n|=2$ for all $n\in \mathbb{N}$,
$$M_p^p(n,C_\phi f) = \frac{1}{|D_n|}\sum\limits_{|v|=n} |f(\phi(v)|^p \leq \frac{2\|f\|^p+2\|f\|^p}{2} =2\|f\|^p,
$$
showing that $\|C_\phi(f)\|^p\leq 2\|f\|^p$ and the result follows.
\epf

\bthm%\label{}
If $T$ is a $(q+1)$-homogeneous tree with $q\geq 2$ such that
\be\label{suff condi}
\sup\limits_{n\in \mathbb{N}} \left(\sum\limits_{|v|=n} q^{|\phi(v)|-n}\right)<\infty,
\ee
then $C_\phi$ is bounded on $\mathbb{T}_p$,  $1\leq p<\infty$.
\ethm
\bpf
For $n\in \mathbb{N}$, $w\in T$ and $f\in \mathbb{T}_p$,   by definition and Lemma \Ref{lem:bound} on growth estimate, we have
\beqq
M_p^p(n,C_\phi f) %&=& \frac{1}{|D_n|}\sum\limits_{|v|=n} |f(\phi(v)|^p\\
&\leq & \frac{1}{|D_n|}\sum\limits_{|v|=n} (q+1)q^{|\phi(v)|-1} \|f\|^p = \sum\limits_{|v|=n} q^{|\phi(v)|-n} \|f\|^p.
\eeqq
Moreover, $M_p^p(0,C_\phi f)=|f(\phi(\textsl{o}))|^p\leq (q+1)q^{|\phi(\textsl{o})|-1} \|f\|^p$
and thus,
$$\|C_\phi f\|^p \leq \max \left\{ (q+1)q^{|\phi(\textsl{o})|-1},~~
\sup\limits_{n\in \mathbb{N}} (\sum\limits_{|v|=n} q^{|\phi(v)|-n}) \right \} \|f\|^p
$$
showing that $C_\phi$ is bounded on $\mathbb{T}_p$.
\epf

\bthm\label{normTp}
Let $T$ be a  $(q+1)$-homogeneous tree and  $1\leq p<\infty$. If $C_\phi$ is bounded on $\mathbb{T}_{p}$,   then
$$\sup\limits_{w\in T} \sup\limits_{n\in \mathbb{N}_0} \left \{ \frac{W(w)}{|D_n|} N_{\phi}(n,w) \right \} \leq \|C_\phi\|^p.
$$
\ethm
\bpf
For each $w\in T$, define $f_w= \left \{ W(w) \chi_{\{w\}} \right \}^{\frac{1}{p}},$ where $W$ is defined in \eqref{MP2-eq1}.
It is easy to verify that for every $w\in T$,
 $M_p(n,f_w)=1$ when $n=|w|$ and $0$ otherwise. This observation gives that $\|f_w\|=1$ for all $w\in T$. Now, for each
 fixed $w\in T$, we have for $n\in \mathbb{N}_0$,
\beqq
M_p^p(n,C_\phi f_w)
&=& \frac{1}{|D_n|}\sum\limits_{|v|=n} W(w) \chi_{\{w\}}(\phi(v))\\
&=& \frac{1}{|D_n|}\sum_{\substack{|v|=n \\ \phi(v)=w }} W(w) = \frac{W(w)}{|D_n|} N_{\phi}(n,w)
\eeqq
which yields that
$$\|C_\phi f_w \|^p = \sup\limits_{n\in \mathbb{N}_0} \left \{ \frac{W(w)}{|D_n|} N_{\phi}(n,w) \right \}.
$$
Consequently,
$$\|C_\phi\|^p = \sup\limits_{\|f\|=1} \|C_\phi(f)\|^p \geq \sup\limits_{w\in T} \|C_\phi(f_w)\|^p
 = \sup\limits_{w\in T} \sup\limits_{n\in \mathbb{N}_0} \left \{ \frac{W(w)}{|D_n|} N_{\phi}(n,w) \right \}
$$
and the desired conclusion follows.
\epf

\bcor%\label{}
If $C_\phi$ is bounded on $\mathbb{T}_{p}$, then
$$\sup \left \{q^{|w|-n} N_{\phi}(n,w): w\in T\setminus \{\textsl{o}\}, n\in \mathbb{N}\right \}
$$ is finite.
\ecor
\bpf
For $w\in T\setminus \{\textsl{o}\}$ and $n\in \mathbb{N}$, we note that
$$\frac{W(w)}{|D_n|}=q^{|w|-n}.
$$
The desired result follows by Theorem \ref{normTp}.
\epf

\bcor
If $\phi$ fixes the root, namely, $\phi(\textsl{o})=\textsl{o}$, then $\|C_\phi\|\geq 1$.
\ecor
\bpf
Let $f$ be the characteristic function on the root $\textsl{o}$. Clearly, $\|f\|=1$ and
$M_p(0,C_\phi f)=|f(\phi(\textsl{o}))|=|f(\textsl{o})|=1$.
We see that
$$\|C_\phi\|=\sup\limits_{\|g\|=1} \|C_\phi(g)\|\geq \|C_\phi(f)\|\geq M_p(0,C_\phi f)=1
$$
and the result follows.
\epf

\bcor\label{Cor:norm-bound}
If $\phi$ does not fix the root, i.e. $\phi(\textsl{o})\neq \textsl{o}$, then
$$\|C_\phi\|^p\geq (q+1) q^{|\phi(\textsl{o})|-1}.
$$
\ecor
\bpf
Let $w=\phi(\textsl{o})$ and, as before,  consider $f_w= \left \{ W(w) \chi_{\{w\}} \right \}^{\frac{1}{p}} $. Now, we observe that
$$\|f_w\|=1~\mbox{ and }~M_p^p(0,C_\phi f_w)=|f_w(\phi(\textsl{o}))|^p= (q+1)q^{|w|-1}
$$
which shows that $\|C_\phi(f_w)\|^p\geq (q+1)q^{|w|-1}$
and the proof follows.
\epf

\bcor\label{Cor:norm}
If $T$ is a $2$-homogeneous tree and $\phi(\textsl{o})\neq \textsl{o}$, then $\|C_\phi\|^p=2$.
\ecor
\bpf
Setting $q=1$ in Corollary \ref{Cor:norm-bound} and Theorem \ref{MP2-th1} gives $\|C_\phi\|^p\geq  2$ and  $\|C_\phi\|^p\leq 2$,
respectively.
\epf

A self map $\phi$ of $T$ is called  an automorphism of $T$, denoted as $\phi \in \Aut (T)$, if $\phi$ is bijective and any
two vertices $v$, $w$ are neighbours $(v \sim w)$ if and only if $\phi(v) \sim \phi(w)$. Now we will compute the norm of the
composition operator $C_\phi$ when the inducing symbol $\phi$ is an automorphism of $T$.

\bthm
Let $T$ be a $(q+1)$-homogeneous tree and consider $C_\phi$ on $\mathbb{T}_{p}$, where $1\leq p<\infty$ and $\phi \in \Aut (T)$. Then we have
\begin{enumerate}
\item[{\rm (i)}]  $\|C_\phi\| = 1$ if
$\phi(\textsl{o})= \textsl{o}$
\item[{\rm (ii)}] $\|C_\phi\|^p = (q+1) q^{|\phi(\textsl{o})|-1}$ if $\phi(\textsl{o})\neq \textsl{o}.$
\end{enumerate}
In particular, every $\phi \in \Aut (T)$ induces bounded composition operator $C_\phi$ on $\mathbb{T}_{p}$.
\ethm
\bpf
Let $D_n=\{v\in T: |v|=n\}$ and consider the case $\phi(\textsl{o})= \textsl{o}$. Then, for each $n$, $\phi$ is a bijective
map from $D_n$ to $D_n$ (since $\phi \in \Aut (T)$ and $\phi(\textsl{o})= \textsl{o}$). For $n\in \mathbb{N}_0$ and
$f\in \mathbb{T}_{p}$, we thus have
\beqq
M_p^p(n,C_\phi f)
&=& \frac{1}{|D_n|}\sum\limits_{|\phi(v)|=n} |f(\phi(v))|^p = M_p^p(n,f).
\eeqq
Taking supremum on both sides, we get $\|C_\phi(f)\|= \|f\|$ which proves the first part.

Next, we consider the case $\phi(\textsl{o})\neq \textsl{o}$. The result is obviously true for $q=1$, by Corollary \ref{Cor:norm}.
Thus, it suffices to prove the theorem for $(q+1)$-homogeneous tree with $q\geq 2$. Let $k=|\phi(\textsl{o})|$.
Since $\phi \in \Aut (T)$, is easy to see that
\begin{center}
\begin{tabular}{|c|c|c|c|c}
\hline
\bf Domain  & \bf ~ Range of $\phi$ contained in& ~~~~ \bf Number of circles\\
\hline
 ~$D_{0}$ & ~$D_{k}$ & $1$ \\
\hline
~$D_{m}$ & ~$D_{k+m}, D_{k+m-2}, \cdots, D_{k-m}$ & $m+1$\\
$(1\leq m\leq k-1)$ & &\\
\hline
~$D_{k}$ & ~$D_{2k}, D_{2k-2}, \cdots ,D_{2},D_{0}$ & $k+1$ \\
\hline
~$D_{k+m+1}$ $(m\geq 0)$ & ~$D_{2k+m+1},D_{2k+m-1}, \cdots ,  D_{2m+1}$ & $k+1$ \\
\hline
\end{tabular}
\end{center}
$$M_p^p(0,C_\phi f)=|f(\phi(0))|^p\leq (q+1) q^{k-1}\|f\|^p.
$$
For the remaining part of the proof, we need to deal with the cases $n=m$ $(1\leq m\leq k-1)$, $n=k$, and  $n\geq k+1$
separately. We begin with
\vspace{8pt}

\noindent{$\ds M_p^p(m,C_\phi f)$}
\beqq
 &=& \frac{1}{(q+1)q^{m-1}}\sum\limits_{|v|=m} |f(\phi(v)|^p\\
& \leq & \frac{1}{(q+1)q^{m-1}} \left \{\sum\limits_{|v|=k+m} |f(v)|^p+\sum\limits_{|v|=k+m-2} |f(v)|^p+
\cdots +\sum\limits_{|v|=k-m} |f(v)|^p \right \} \\
&\leq & \frac{1}{(q+1)q^{m-1}} \left \{ (q+1) q^{k+m-1}+ (q+1) q^{k+m-3} + \cdots +(q+1) q^{k-m-1} \right \}\|f\|^p\\
&=& \left \{q^{k} +  q^{k-2} + \cdots + q^{k-2m}  \right \} \|f\|^p\\
&\leq & (q+1) q^{k-1} \|f\|^p
\eeqq
showing that $M_p^p(n,C_\phi f) \leq  (q+1) q^{k-1} \|f\|^p$ for $n=1,2, \cdots ,k-1.$ Next,  for $n=k$, we find that

\vspace{8pt}

\noindent{$\ds M_p^p(k,C_\phi f)$}
\beqq
~~~~~~~~  &=& \frac{1}{(q+1)q^{k-1}}\sum\limits_{|v|=k} |f(\phi(v)|^p\\
& \leq & \frac{1}{(q+1)q^{k-1}} \left \{\sum\limits_{|v|=2k} |f(v)|^p+\sum\limits_{|v|=2k-2} |f(v)|^p+\cdots +
\sum\limits_{|v|=2} |f(v)|^p+|f(\textsl{o})|^p \right \} \\
&\leq & \frac{1}{(q+1)q^{k-1}} \left \{ (q+1) q^{2k-1}+ (q+1) q^{2k-3} + \cdots +(q+1) q + 1 \right \}\|f\|^p\\
&=& \left \{q^{k} +  q^{k-2} + \cdots + q^{2-k} + \frac{1}{(q+1) q^{k-1}} \right \} \|f\|^p\\
&\leq& \left \{q^{k} +  q^{k-2} + \cdots + q^{2-k} +  q^{1-k} \right \} \|f\|^p\\
&\leq & (q+1) q^{k-1} \|f\|^p.
\eeqq
Finally, for each $m \in \mathbb{N}_{0}$,

\vspace{8pt}

\noindent{$\ds M_p^p(m+k+1,C_\phi f) $}
\beqq
  &=& \frac{1}{(q+1)q^{m+k}}\sum\limits_{|v|=m+k+1} |f(\phi(v)|^p\\
& \leq & \frac{1}{(q+1)q^{m+k}} \left \{\sum\limits_{|v|=m+2k+1} |f(v)|^p+\sum\limits_{|v|=m+2k-1} |f(v)|^p+\cdots
+\sum\limits_{|v|=2m+1} |f(v)|^p \right \} \\
&\leq & \frac{1}{(q+1)q^{m+k}} \left \{ (q+1) q^{m+2k}+ (q+1) q^{m+2k-2} + \cdots +(q+1) q^m \right \}\|f\|^p\\
&=& \left \{q^{k} +  q^{k-2} + \cdots + q^{-k}  \right \} \|f\|^p\\
&\leq & (q+1) q^{k-1} \|f\|^p.
\eeqq
The above discussion implies that
$$M_p^p(n,C_\phi f) \leq  (q+1) q^{k-1} \|f\|^p ~\mbox{ for all $n \in \mathbb{N}_{0}$}
$$
and thus, $\|C_\phi\|^p \leq (q+1) q^{|\phi(\textsl{o})|-1}$.
Other way inequality follows from Corollary \ref{Cor:norm-bound} and the proof is complete.
\epf
\section{Compact Composition Operators} \label{MP2Sec5}

A bounded linear operator $A$ from a normed linear space $X$ to a normed linear space $Y$ is said to be a
\textit{compact operator} if the image of closed unit ball $ \{Ax: \|x\|_X\leq1\}$ has compact closure in $Y$.

In the classical case, for an analytic self map $\phi$ of $\D$, the following statements are
equivalent (see \cite[Section 2.7 and Compactness Theorem, Chapter 10]{Shapiro:Book}):
\bee
\item[(a)] $C_\phi$ is compact on $H^p$  for $1\leq p<\infty$.
\item[(b)]  $C_\phi$ is compact on $H^2$.
\item[(c)] $\ds\lim\limits_{|w|\rightarrow1^-} \ds{\frac{N_{\phi}(w)}{\log\frac{1}{|w|}}}=0,$
where $N_{\phi}$ is the Nevanlinna counting function of $\phi$.
\eee
Also, $C_\phi$ is compact on $H^{\infty}$ if and only if $\sup\{|\phi(z)|:z\in\D\}<1$
(see \cite[Problem 10, Chapter 2]{Shapiro:Book}).

For the discrete setting,  we now consider the compactness of composition operators on $\mathbb{T}_{p}$ spaces.
A self map $\phi$ of $(q+1)$-homogeneous tree $T$ is called a bounded map if $\sup \{|\phi(v)| :v\in T\}$ is finite.

\bthm\label{lem:cmpt}
Every bounded self map $\phi$ of $T$ induces compact composition operator on $\mathbb{T}_{p}$ for $1\leq p\leq\infty$.
\ethm
\bpf
Suppose $\phi$ is a bounded self map of a $(q+1)$-homogeneous tree $T$. Then ${\rm Range}\,(\phi)$ is finite set,
say, ${\rm Range}\,(\phi)=\{v_1,v_2,\ldots,v_k\}$. For each $1\leq i\leq k$, denote by $E_i$ for the pre-image of $v_i$
under $\phi$. If $\phi(v)=v_i$, then $f\circ\phi(v)=f(v_i)$ so that
$$f\circ\phi = f(v_1)\chi_{E_1}+f(v_2)\chi_{E_2}+\cdots+f(v_k)\chi_{E_k}
$$
and ${\rm Range}\,(C_\phi)={\rm span}\,\{\chi_{E_1},\chi_{E_2},\ldots,\chi_{E_k}\}$. Thus, $C_\phi$ is a finite rank
operator and hence it is compact.
% (see \cite[Theorem 4.18(a)]{Rudin:Book}).
\epf

\bthm\label{MP2-th6}
If  $\phi$ is a self map of $(q+1)$-homogeneous tree $T$, then the following are equivalent:
\begin{enumerate}
\item[{\rm (a)}] $C_\phi$ is compact on $\mathbb{T}_{p}$ for $1\leq p\leq\infty$.
\item[{\rm (b)}] $\|C_\phi f_n\| \rightarrow 0$ as $n\rightarrow \infty$ whenever bounded sequence of functions
$\{f_n\}$ that converges to $0$ pointwise.
\end{enumerate}
\ethm
\bpf
(a)~$\Ra$~(b):
Assume that $C_\phi$ is compact on $\mathbb{T}_{p}$ and $\{f_n\}$ is a bounded sequence in $\mathbb{T}_{p}$ that
converges to $0$ pointwise. Suppose on the contrary that $\|C_\phi(f_n)\| \not\rightarrow 0$ as $n \rightarrow \infty$.
Then there exists a subsequence $\{f_{n_j}\}$ and an $\epsilon >0$ such that
$\|C_\phi(f_{n_j})\| \geq \epsilon$ for all $j$. Denote $\{f_{n_j}\}$ by $\{g_j\}$.
Since $C_\phi$ is compact, there is a subsequence $\{g_{j_k}\}$ of $\{g_{j}\}$ such that $ \{C_\phi (g_{j_k})\}$ converges
to some function, say, $g$.
It follows that $\{ C_\phi (g_{j_k})\}$ converges to $g$ pointwise
% . Since ${g_j}$ converges pointwise to $0$ forces that
and $g \equiv 0$ implying that $\{C_\phi (g_{j_k})\}$ converges to $0$ which is a contradiction to $\|C_\phi(g_j)\| \geq \epsilon$ for all $j$.
Hence, $\|C_\phi(f_n)\| \rightarrow 0$ as $n\rightarrow\infty$.

(b)~$\Ra$~(a): Conversely, suppose that case (b) holds. First let us consider the case $1\leq p<\infty$. Let $\{g_n\}$ be a
sequence in unit ball of $\mathbb{T}_{p}$. By Lemma \Ref{lem:bound}, for each $v\in T$, the sequence $\{g_n(v)\}$ is bounded.
By the diagonalization process, there is a subsequence
$\{g_{nn}\}$ of $\{g_n\}$ that converges pointwise to g (say).
We see that, for each $m \in \mathbb{N}_{0}$,
\beqq
M_p^p(m,g)
&=& \lim\limits_{n\rightarrow\infty} \frac{1}{(q+1)q^{m-1}}\sum\limits_{|v|=m} |g_{nn}(v)|^p
\leq  \limsup \|g_{nn}\|^p \leq 1
\eeqq
showing that $g\in \mathbb{T}_{p}$ with $\|g\|\leq 1$. Consequently, if $f_n=g_{nn}-g$, then $\{f_n\}$ converges
to $0$ pointwise and $\|f_n\|\leq 2$. By the assumption (b), $\|C_\phi f_n\| \rightarrow 0$ as $n\rightarrow \infty$
and thus, $\{C_\phi(g_{nn})\}$ converges to $C_\phi(g)$.  Hence $C_\phi$ is compact on $\mathbb{T}_{p}$.

 The proof for the case $p=\infty$ is similar to the above.
\epf

\brem {\rm
Since edge counting metric on $T$ induces discrete topology, compact sets are only sets having finitely many
elements. Thus convergence uniformly on compact subsets of $T$ is equivalent to pointwise convergence. In view of this
remark, Theorem \ref{MP2-th6} is a discrete analog of weak
convergence theorem (see \cite[section 2.4, p.~29]{Shapiro:Book}) in the classical case.
}
\erem
\bcor%\label{}
Let $\phi$ be a self map of $T$. Then $C_\phi$ is compact on $\mathbb{T}_{\infty}$ if and only if $\phi$ is a
bounded self map of $T$.
\ecor
\bpf
 If $\phi$ is a bounded self map of $T$, then $C_\phi$ is compact, by Theorem \ref{lem:cmpt}. Conversely,
 suppose $\phi$ is not a bounded map. Then, there exists a sequence of vertices $\{v_k\}$ of $T$ such that
$\phi(v_k)=w_k$ and $|w_k|\rightarrow\infty$ as $k\rightarrow\infty$. Take $f_k=\chi_{\{w_k\}}$ for each $k\in\N$.
Then, $\|f_k\|_\infty=1$ for each $k$ and $\{f_k\}$ converges to $0$ pointwise. Since $C_\phi$ is compact,
$\|C_\phi(f_k)\|_\infty \rightarrow 0$ as $k\rightarrow \infty$, by Theorem \ref{MP2-th6}. This is not possible, because $\|C_\phi(f_k)\|_\infty=1$
for each $k\in\N$, which can be observed from the definition of $f_k$. Hence $\phi$ should be a bounded map.
\epf
\bcor\label{MP2-cor5}
Let $T$ be a  $(q+1)$-homogeneous tree and $1\leq p<\infty$. If $C_\phi$ is compact on $\mathbb{T}_{p}$, then
$$\sup\limits_{n\in \mathbb{N}_0} \left \{ q^{|w|-n} N_{\phi}(n,w) \right \} \rightarrow 0 ~\mbox{ as}~ |w|\rightarrow\infty.
$$
\ecor\bpf
As in the earlier situations, for each $w\in T\setminus \{\textsl{o}\}$, we let $f_w= \left \{ W(w) \chi_{\{w\}} \right \}^{\frac{1}{p}}$.
% where $\chi_{\{v\}}$ is characteristic function on $\{v\}$.
Then, $\|f_w\|=1$ for all $w$ and, since $f_w(v)=0$ whenever $|w|> n=|v|$,  $\{f_w\}$ converges
to $0$ ponitwise. Since $C_\phi$ is compact, $\|C_\phi (f_w)\| \rightarrow 0 ~\mbox{ as }~ |w|\rightarrow\infty$.
However, we have already shown that
$$\|C_\phi f_w \|^p = \sup\limits_{n\in \mathbb{N}_0} \left \{ \frac{W(w)}{|D_n|} N_{\phi}(n,w) \right \}= \sup\limits_{n\in \mathbb{N}_0} \left \{ q^{|w|-n} N_{\phi}(n,w) \right \}
$$
and the desired conclusion follows.
\epf

\brem
{\rm
For $2$-homogeneous trees,  Corollary \ref{MP2-cor5} takes a simpler form:  If $C_\phi$ is compact  on $\mathbb{T}_{p}$, then
$$ \sup\limits_{n\in \mathbb{N}_0} \left \{  N_{\phi}(n,w) \right \} \rightarrow 0 ~\mbox{ as }~  |w|\rightarrow\infty .
$$
This remark is helpful in the proof of Corollary \ref{MP2-th7}.
}
\erem

\bcor\label{MP2-cor6}
If $C_\phi$ is compact on $\mathbb{T}_{p}$, then $|v|-|\phi(v)|\rightarrow \infty$ as $|v|\rightarrow\infty$.
\ecor
\bpf
We will prove this result by contradiction. Suppose that $|v|-|\phi(v)| \not\rightarrow \infty$ as $|v|\rightarrow\infty$.
Then there exists a sequence of vertices $\{v_{k}\}$ and an $M >0$ such that $ |v_k| -|\phi(v_k)| \leq M$ for
all $k$ which implies  that
$|\phi(v_k)|\rightarrow\infty$ as $k\rightarrow\infty$. Since $N_{\phi}(|v_{k}|,\phi(v_{k}))\geq 1$ for all $k$,
where  $N_{\phi}(n,w)$ is defined as
in Section \ref{MP2Sec4}, we obtain that $N_{\phi}(|v_{k}|,\phi(v_{k})) q^{|\phi(v_{k})|-|v_{k}|}\geq q^{-M}$
which yields that
$$\sup\limits_{n\in \mathbb{N}_0} \left \{ N_{\phi}(n,\phi(v_{k})) q^{|\phi(v_{k})|-n} \right \}\geq q^{-M} ~\mbox{ for all $k$}
$$
and thus,
$$\sup\limits_{n\in \mathbb{N}_0} \left \{ q^{|w|-n} N_{\phi}(n,w) \right \} \not\rightarrow 0 ~\mbox{ as }~  |w|\rightarrow\infty
$$
which gives that $C_\phi$ is not compact, by Corollary \ref{MP2-cor5}. This contradiction completes the proof.
\epf

\bcor\label{MP2-th7}
Let $T$ be a $2$-homogeneous tree. Then $C_\phi$ is compact on $\mathbb{T}_{p}$ if and only if $\phi$ is a bounded self map of $T$.
\ecor
\bpf
Since every bounded self map $\phi$ of $T$ induces compact composition operator on $\mathbb{T}_{p}$, one way implication is hold.
For the proof of the converse part, we suppose that $\phi$ is not bounded. Then the range contains an infinite set, say $\{w_1,w_2,\ldots\}$.
For each $k$, choose $v_k \in T$ such that $\phi(v_k)=w_k$. This gives  $N_{\phi}(|v_{k}|,w_{k}) \geq 1$ and thus
$\sup\limits_{n\in \mathbb{N}_0} N_{\phi}(n,w_{k}) \geq 1 $ for all $k$. It follows that
$\sup\limits_{n\in \mathbb{N}_0} \left \{  N_{\phi}(n,w) \right \} \not\rightarrow 0 ~\mbox{ as }~  |w|\rightarrow\infty$
and hence, $C_\phi$ cannot be compact.
\epf
\brem {\rm
It is worth to recall from \cite[Chapter 3, p.~37]{Shapiro:Book} that if a ``big-oh" condition describes a
class of bounded operators, then the corresponding ``little-oh" condition picks out the subclass of compact operators".
We have already shown that if $\sum _{|v|=n} q^{|\phi(v)|}=O(q^n)$ then $C_\phi$ is bounded on $\mathbb{T}_{p}$.
So it is natural to ask whether $\sum _{|v|=n} q^{|\phi(v)|}=o(q^n)$ guarantees the compactness of $C_\phi$  on $\mathbb{T}_{p}$.
Indeed, the answer is yes. Clearly the later observation is not useful because no self map $\phi$ of $T$ satisfies this condition.
This is because $\sum _{|v|=n} q^{|\phi(v)|}\geq \sum _{|v|=n} q^{0} =(q+1) q^{n-1}$ and thus,
$\sum_{|v|=n} q^{|\phi(v)|}=o(q^n)$ is cannot be possible.
}
\erem

\section{Examples}\label{MP2Sec6}
\beg
{\rm
 For each $n\in \mathbb{N}_0$, choose the vertex $v_n$ such that $v_n\in D_n=\{v\in T: |v|=n\}$. Define $\phi_1(v)=v_n$ if $|v|=n$.
% Now we consider the composition operator induced by $\phi_1$.
Now, for $q\neq1$, consider the function $f$ defined by
$$f(v)= \left \{\begin{array}{cl}
\ds \{(q+1)q^{n-1}\}^{\frac{1}{p}} & \mbox{ if $v=v_n$ for some $n\in \N$}\\
0 &\mbox{ elsewhere.}
\end{array} \right .
$$
Then $f \in \mathbb{T}_{p}$, $\|f\|=1$ and for each $m \in \mathbb{N}$, we see that
\beqq
M_p^p(m,C_{\phi_1} f)
%&=& \frac{1}{(q+1)q^{m-1}}\sum\limits_{|v|=m} |f(\phi_1(v))|^p\\
&=& \frac{1}{(q+1)q^{m-1}}\sum\limits_{|v|=m} |f(v_m)|^p = (q+1)q^{m-1}
\eeqq
showing that $\|C_{\phi_1}f\| = \sup\limits_{m\in \mathbb{N}_{0}} M_{p}(m,C_{\phi_1}f)$ which is not finite for $q\geq 2$.
This example shows that there are self maps of $T$ which do not induce bounded composition operator on
$\mathbb{T}_{p}$ unlike the case of Hardy spaces on the unit disk.
}
\eeg

\beg {\rm
Consider the following self map $\phi_2$ of $T$ defined by
$$\phi_2(v)= \left \{\begin{array}{cl}
\textsl{o} & \mbox{ if $v=\textsl{o}$}\\
v^- &\mbox{ otherwise.}
\end{array} \right .
$$
where $v^-$ denotes the parent of $v$. Then it follow easily that
$$ M_p^p(0,C_{\phi_2} f) = M_p^p(0,f)~\mbox{ and }~ M_p^p(1,C_{\phi_2} f)  = \frac{1}{(q+1)}\sum\limits_{|v|=1} |f(\textsl{o})|^p  = M_p^p(0,f).
$$
Finally, for $n\geq2$, we have
\beqq
M_p^p(n,C_{\phi_2} f)  &=& \frac{1}{(q+1)q^{n-1}}\sum\limits_{|v|=n} |f(v^-)|^p\\
&=& \frac{q}{(q+1)q^{n-1}}\sum\limits_{|w|=n-1} |f(w)|^p\\
&=& M_p^p(n-1,f)
\eeqq
and thus,
$$\|C_{\phi_2}f\| = \sup\limits_{m\in \mathbb{N}_{0}} M_{p}(m,C_{\phi_2}f) = \sup\limits_{m\in \mathbb{N}_{0}} M_{p}(m,f) = \|f\|
$$
showing that  $C_{\phi_2}$ is bounded on $\mathbb{T}_{p}$. On the other hand, since $|v|-|\phi_2(v)|=1$ for all $|v|\geq1$, we have
$|v|-|\phi_2(v)|\not\rightarrow \infty$ as $|v|\rightarrow\infty$. Hence $C_{\phi_2}$  is not compact, by Corollary \ref{MP2-cor6}.
This is an example of bounded composition operator on $\mathbb{T}_{p}$ which is not compact.
}
\eeg

\brem {\rm
Let $\phi_3$ be a map on $T$ such that $\phi_3$ maps every vertex into any one of its child. Then, as in the case
of $C_{\phi_2}$,  it is easy to see that $C_{\phi_3}$ is bounded but not compact.
Moreover,  it can be seen that, for each $n\in \mathbb{N}$, $(C_{\phi_2})^n$ and $(C_{\phi_3})^n$ are also bounded
but not compact.
}
\erem

\beg {\rm
For each $n\in \mathbb{N}_0$, choose a vertex $v_n$ such that $|v_n|=n$. Define a self map $\phi_4$ by
$$\phi_4(v)= \left \{\begin{array}{cl}
v_k & \mbox{ if $v=v_{2k}$ for some $k\in \N$}\\
\textsl{o} &\mbox{ otherwise.}
\end{array}\right .
$$
Then we obtain that
$$ M_p^p(0,C_{\phi_4} f) = |f(\phi_4(\textsl{o}))|^p = |f(\textsl{o})|^p = M_p^p(0,f).
$$
Next, for an odd natural number $n$, we see that
%\beqq
%M_p^p(n,C_{\phi_4} f) &=& \frac{1}{(q+1)q^{n-1}}\sum\limits_{|v|=n} |f(\phi_4(v))|^p\\
%&=& \frac{1}{(q+1)q^{n-1}}\sum\limits_{|v|=n} |f(\textsl{o})|^p\\
%&=& |f(\textsl{o})|^p
%\eeqq
$$M_p^p(n,C_{\phi_4} f)  = \frac{1}{(q+1)q^{n-1}}\sum\limits_{|v|=n} |f(\textsl{o})|^p = |f(\textsl{o})|^p.
$$
Finally, for an even natural number, say $n=2k$, for some $k\in \mathbb{N}$, we find that
\beqq
M_p^p(n,C_{\phi_4} f) % &=& \frac{1}{(q+1)q^{n-1}}\sum\limits_{|v|=n} |f(\phi_4(v))|^p\\
&=& \frac{1}{(q+1)q^{n-1}}\left\{\sum_{\substack{|v|=n \\ v\neq v_{2k} }} |f(\phi_4(v))|^p + |f(\phi_4(v_{2k}))|^p \right\}\\
&\leq& |f(\textsl{o})|^p + \frac{|f(v_{k})|^p }{(q+1)q^{n-1}}.\\
\eeqq
Thus, by Lemma \Ref{lem:bound}, we have
$$\|C_{\phi_4}f\|^p \leq |f(\textsl{o})|^p + \sup\limits_{k\in \mathbb{N}} \left\{\frac{f(v_{k})|^p}{(q+1)q^{2k-1}} \right \} \leq 2 \|f\|^p
$$
which shows that  $C_{\phi_4}$ is bounded  on $\mathbb{T}_{p}$.

Suppose now that $T$ is a $(q+1)$-homogeneous tree with $q\geq 2$. Let $\{f_n\}$ be a
sequence in the unit ball of $\mathbb{T}_{p}$ which converges to $0$ pointwise. Note that
$\left\{\frac{f(v_{k})|^p}{(q+1)q^{2k-1}} \right \} \leq \frac{1}{q^k}$,
by Lemma \Ref{lem:bound}. We now claim that $\|C_{\phi_4}f_n\|^p \rightarrow 0$ as $n\rightarrow\infty$.

Let $\epsilon > 0$ be given. Then there exists a natural number $N_1$ such that $q^{-k} < \epsilon/2$ for all $k\geq N_1$.
Consider the set $S=\{v_1,v_2,\ldots,v_{N_1}\}$. Since $\{f_n\}$ converges to $0$ pointwise, we can choose a natural
 number $N>N_1$ such that $|f_n(\textsl{o})|^p <\epsilon/{2}$
and $|f_n(v)|^p <\epsilon/{2}$ for all $v\in S$ and for all $n\geq N$. Thus,
$$\|C_{\phi_4}f_n\|^p \leq |f_n(\textsl{o})|^p + \sup \left\{ \frac{\epsilon}{2}, \frac{1}{q^{N_1}},\frac{1}{q^{N_{1}+1}},\cdots \right \} \leq \epsilon ~\mbox{ for all }~ n\geq N
$$
which gives that $\|C_{\phi_4}f_n\|^p \rightarrow 0$ as $n\rightarrow\infty$ and hence $C_{\phi_4}$ is compact on $\mathbb{T}_{p}$.
This example shows that there are unbounded self maps of $T$ which induces compact composition operators on
$\mathbb{T}_{p}$ for the case of $(q+1)$-homogeneous trees with $q\geq 2$.
}
\eeg

We conclude the paper with a comparison in the case of $q=1$ and $q>1$. For $2$-homogeneous trees, circle of radius
$n$ has $2$ vertices for all $n\in \mathbb{N}$ whereas in the case of $(q+1)$-homogeneous trees with $q\geq 2$,
circle of radius $n$ has $(q+1) q^{n-1}$ vertices. Due to this basic fact, we can expect a difference in operator
theoretic point of view. The following table explains how composition operators on $\mathbb{T}_{p}$ for $1\leq p<\infty$
differ in both the cases.\\

\begin{center}
\begin{tabular}{|l|l|}
\hline
\bf $2$-homogeneous tree  & \bf $(q+1)$-homogeneous tree with $q\geq 2$\\
\hline
Every self map $\phi$ of $T$ induces & There are self maps of $T$ which induces\\
bounded composition & unbounded composition\\
operators on $\mathbb{T}_{p}$ &   operators on $\mathbb{T}_{p}$.\\
\hline
Only bounded self map of $T$   &  There are unbounded self maps of $T$ \\
induces compact composition  &  which induces compact operators  \\
operators on $\mathbb{T}_{p}$ &  on $\mathbb{T}_{p}$.\\
\hline
\end{tabular}
\end{center}

\subsection*{\bf Acknowledgement}
The authors thank the referee for many useful comments.
The first author thanks the Council of Scientific and Industrial Research (CSIR), India,
for providing financial support in the form of a SPM Fellowship to carry out this research.
The second  author is currently on leave from the IIT Madras.

\end{document}